\documentclass[review,onefignum,onetabnum]{siamart171218}

\usepackage{amsfonts}
\usepackage{graphicx}
\usepackage{epstopdf}
\usepackage{algorithmic}
\ifpdf
  \DeclareGraphicsExtensions{.eps,.pdf,.png,.jpg}
\else
  \DeclareGraphicsExtensions{.eps}
\fi

\usepackage{enumitem}
\setlist[enumerate]{leftmargin=.5in}
\setlist[itemize]{leftmargin=.5in}

\newsiamremark{remark}{Remark}
\newsiamremark{hypothesis}{Hypothesis}
\crefname{hypothesis}{Hypothesis}{Hypotheses}
\newsiamthm{claim}{Claim}
\newtheorem{example}[theorem]{Example}

\headers{Regularization of A Coupled System with White Gaussian Noise}{N. H. Tuan, V. A. Khoa, P. T. K. Van, and V. V. Au}

\title{An Improved Quasi-Reversibility Method for A Terminal-Boundary Value Multi-species Model with White Gaussian Noise
\thanks{Submitted to the editors DATE.
\funding{The second author was funded by US Army Research Laboratory and US Army Research Office grant W911NF-19-1-0044. Also, the work of the second author was also partly supported by the Research Foundation-Flanders (FWO) under the project named ``Approximations for forward and inverse reaction-diffusion problems related to cancer models''.}}}

\author{Nguyen Huy Tuan\thanks{Applied Analysis Research Group, Faculty of Mathematics and Statistics, Ton Duc Thang University, Ho Chi Minh City, Vietnam 
  (\email{nguyenhuytuan@tdtu.edu.vn}).}
\and Vo Anh Khoa\thanks{Corresponding author. Department of Mathematics and Statistics, University of North Carolina at Charlotte, Charlotte, North Carolina 28223, USA, and Faculty of Sciences, Hasselt University, Campus Diepenbeek, Agoralaan Building D, BE3590 Diepenbeek, Belgium. 
	(\email{vakhoa.hcmus@gmail.com,anhkhoa.vo@uncc.edu}).}
\and Phan Thi Khanh Van\thanks{Faculty of Applied Science, Ho Chi Minh City University of Technology, Ho Chi Minh City, Vietnam, and Faculty of Mathematics and Computer Science, University of Science, Vietnam National University, 227 Nguyen Van Cu, District 5, Ho Chi Minh City, Vietnam
		(\email{khanhvanphan@hcmut.edu.vn}).}
\and Vo Van Au\thanks{Institute of Fundamental and Applied Sciences, Duy Tan University, Ho Chi Minh City 700000, Vietnam, and  Faculty of Natural Sciences, Duy Tan University, Da Nang, 550000, Vietnam. (\email{vovanau@duytan.edu.vn}).}
}

\usepackage{amsopn}

\ifpdf
\hypersetup{
  pdftitle={Regularization of A Coupled System with White Gaussian Noise},
  pdfauthor={N. H. Tuan, V. A. Khoa, P. T. K. Van, and V. V. Au}
}
\fi

\begin{document}

\maketitle

\begin{abstract}
  Upon the recent development of the quasi-reversibility method for terminal value parabolic problems in \cite{Nguyen2019}, it is imperative to investigate the convergence analysis of this regularization method in the stochastic setting. In this paper, we positively unravel this open question by focusing on a coupled system of Dirichlet reaction-diffusion equations with additive white Gaussian noise on the terminal data. In this regard, the approximate problem is designed by adding the so-called perturbing operator to the original problem and by exploiting the Fourier reconstructed terminal data. By this way, Gevrey-type source conditions are included, while we  successfully maintain the logarithmic stability estimate of the corresponding stabilized operator, which is necessary for the error analysis. As the main theme of this work, we prove the error bounds for the concentrations and for the concentration gradients, driven by a large amount of weighted energy-like controls involving the expectation operator. Compared to the classical error bounds in $L^2$ and $H^1$ that we obtained in the previous studies, our analysis here needs a higher smoothness of the true terminal data to ensure their reconstructions from the stochastic fashion. Two numerical examples are provided to corroborate the theoretical results.
\end{abstract}



\section{Introduction}
The main purpose of this paper is to carry out the error analysis of the recently proposed quasi-reversibility (QR) method in a stochastic setting for a class of terminal-boundary value multi-species model. As concluded in \cite{Nguyen2019}, the error bounds we obtained there can be very helpful in the finite element settings due to the variational framework we choose. Nonetheless, it still questions us about if this approach can be modified to somewhat get the convergence analysis in some specific stochastic setting involved in the partial differential equations (PDEs) we are chasing. In this work, we provide a positive answer by taking into account a nonlocal coupled system of nonlinear reaction-diffusion equations with additive white Gaussian noise on the terminal data. Besides, we significantly strengthen the applicability of this general QR framework by introducing a ready-to-use scheme from our computational standpoint.


\subsection{Statement of the problem}
Let $0<T< \infty$ be the final time of observation and $\Omega\subset\mathbb{R}^{d}$ for $\mathbb{N}^{*}\ni d\le 3$ be the domain of interest, which is open, connected, and bounded with a sufficiently smooth boundary. Denoted by $Q_{T}=\Omega \times (0,T)$, this work is devoted to finding $u,v:\overline{Q_{T}}\to \mathbb{R}$ as two solutions of the following evolution system:
\begin{equation}\label{system1}
\begin{cases}
u_{t}-\mathcal{D}_{1}\left(\ell_{0}\left(u\right)\left(t\right)\right)\Delta u=F\left(x,t;u;v\right),\\
v_{t}-\mathcal{D}_{2}\left(\ell_{0}\left(v\right)\left(t\right)\right)\Delta v=G\left(x,t;u;v\right) & \text{for }\left(x,t\right)\in Q_{T}.
\end{cases}
\end{equation}
To complete the terminal-boundary value problem we wish to solve, the system \eqref{system1} is then supplemented with the terminal conditions $u(x,T)=u_{f}(x)$, $v(x,T)=v_{f}(x)$ and the Dirichlet boundary conditions $u = v = 0$ on the boundary $\partial\Omega$. Suppose further that the terminal data are observed in the presence of white noise processes $\xi_{1}$, $\xi_{2}$ in terms of
\begin{align}\label{eq:noise}
u_{f}^{\varepsilon}\left(x\right)=u_{f}\left(x\right)+\varepsilon\xi_{1}\left(x\right),\quad v_{f}^{\varepsilon}\left(x\right)=v_{f}\left(x\right)+\varepsilon\xi_{2}\left(x\right),
\end{align}
where $\varepsilon\in (0,1)$ is used for measuring the amplitude of the noise.

In this scenario, our attempt is to seek the initial values $u(x,0)=u_{0}(x)$ and $v(x,0)=v_{0}(x)$ when we only know the measurements of the terminal data $u_{f}$ and $v_{f}$. We accentuate that this seeking can only be done by regularization since our problem is, in general, exponentially unstable; cf. \cite{Alquier2011,Cavalier2004,Golubev1999,Agapiou2014,Knapik2013,Kaipio2011} for a glimpse of statistical inverse problems. Starting from the source reconstruction model of diffusive, competitive and reactive brain tumor cells (cf. \cite{Jaroudi2018} and references cited therein), solving \eqref{system1} is essential to understand the localization of the tumor source and then to possibly advance the treatment of brain tumors. This biological context conjures up the following images. Let $u$ and $v$ be the normal (healthy) and abnormally growing tissue cells densities in a brain region, respectively. In the perspective of \eqref{system1}, we assume that the movements of each kind of cells are dominantly influenced by the whole population of the corresponding type. Usually, this significant impact brings us to the nonlocal form of diffusion (see \cite{Chipot1997,Almeida2016,Duque2016} for the mathematical background of this diffusion in the forward model), which reads as
\[
\mathcal{D}_{i}\left(\ell_{0}\left(w\right)\left(t\right)\right) = \mathcal{D}_{i}\left(\int_{\Omega}f(x)w(x,t)dx\right),\; \ell_{0}(w)(t):= \int_{\Omega}f(x)w(x,t)dx,
\]
for $i=1,2$, $w\in\mathbb{R}$ and where $f$ is a given sufficiently smooth weight function. Largely motivated by the typical expression, we considered in \cite{Tuan2017} a slight generalization of this diffusion for the sole purpose of extending the applicability of the QR method we are studying.

The next concern that we would point out in the continuous model \eqref{system1} lies in the nonlinearities $F$ and $G$ which can be considered as reaction, death and proliferation rates involved in the active network under scrutiny. In fact, we are able to suppose that at the low-grade gliomas regime, one knows that the brain tumor cells start growing slowly by infiltrating into the healthy brain cells and simultaneously, they misdirect the T cells of the immune system to avoid being beaten. Besides, we can assume that for some vigorous people, those gliomas are still identified at the beginning and then are very slightly removed by some reaction rates from the immune system. Accordingly, this entire incident is well-suited to the consideration of the coupled system \eqref{system1}, provided that the assumptions on these nonlinear terms will cover most of the real-world contexts. Needless to say, one may also detail several well-known chemical models contained in \eqref{system1} such as the two-species Lotka--Volterra competition-diffusion model, the Brusselator coupled system for the autocatalytic reaction and those already mentioned in \cite{Nguyen2019}. Eventually, we remark that even though the evolution system \eqref{system1} can be expressed in a closed-form, it is better to contemplate the whole mathematical treatment for every single term in there.

\subsection{Settings of the statistical terminal data}
Cf. \cite{Alquier2011}, we recall standard assumptions on the stochastic noise we want to address in the model. In the following, $\mathcal{H}$ is a Hilbert space.

\begin{definition}
	The stochastic error is a Hilbert-space process, i.e. a bounded linear operator $\xi:\mathcal{H}\to L^2(\tilde{\Omega},\mathcal{A},P)$ where $(\tilde{\Omega},\mathcal{A},P)$ is a complete probability space and $L^2(\cdot)$ is the space of all square integrable measurable functions.
\end{definition}

This way we are able to define the random variables $\left\langle \xi,g_{j}\right\rangle $ for $j=1,2$ for all $g_{1},g_{2}\in\mathcal{H}$
by definition $\mathbb{E}\left\langle \xi,g_{j}\right\rangle =0$.
Furthermore, we define its covariance $\text{Cov}_{\xi}$ as the bounded
linear operator mapping from $\mathcal{H}$ onto itself such that
$\left\langle \text{Cov}_{\xi}g_{1},g_{2}\right\rangle =\text{Cov}\left(\left\langle \xi,g_{1}\right\rangle ,\left\langle \xi,g_{2}\right\rangle \right)$.

\begin{definition}
	We say that $\xi$ is a \emph{white noise process} in $\mathcal{H}$,
	if $\text{Cov}_{\xi}=I$ and the induced random variables are Gaussian:
	for all $g_{1},g_{2}\in\mathcal{H}$, the random variables $\left\langle \xi,g_{j}\right\rangle $
	have distributions $\mathcal{N}\left(0,\left\Vert g_{j}\right\Vert ^{2}\right)$
	and $\text{Cov}\left(\left\langle \xi,g_{1}\right\rangle ,\left\langle \xi,g_{2}\right\rangle \right)=\left\langle g_{1},g_{2}\right\rangle $.
\end{definition}

Accordingly, assume that the observations \eqref{eq:noise} can only be obtained in a discretized or binned form. This means that we only have vectors of normally distributed random variables $\left\{ u_{f}^{\varepsilon,j}\right\} _{j=\overline{1,n}}$, $\left\{ v_{f}^{\varepsilon,j}\right\} _{j=\overline{1,n}}$ given by
\begin{align}\label{eq:data1}
	u_{f}^{\varepsilon,j} & :=\left\langle u_{f}^{\varepsilon},\phi_{j}\right\rangle =\left\langle u_{f},\phi_{j}\right\rangle +\varepsilon\left\langle \xi_{1},\phi_{j}\right\rangle ,\\\label{eq:data2}
	v_{f}^{\varepsilon,j} & :=\left\langle v_{f}^{\varepsilon},\phi_{j}\right\rangle =\left\langle v_{f},\phi_{j}\right\rangle +\varepsilon\left\langle \xi_{2},\phi_{j}\right\rangle,
\end{align}
where $n\in\mathbb{N}$ is the number of steps of discrete observations and $\phi_j$  is itself taken from the orthonormal basis $\left\{ \phi_{j}\right\} _{j\in\mathbb{N}}$ of $L^2(\Omega)$. Recall that due to the smoothness of $\Omega$, the existence of this basis is guaranteed, where $\phi_j \in H_0^1(\Omega)\cap C^{\infty}(\overline{\Omega})$ solves the basic eigenvalue problem $-\Delta\phi_{j}(x) = \mu_{j}\phi_{j}(x)$ for $x\in\Omega$. Additionally, the Dirichlet eigenvalues  $\left\{ \mu_{j}\right\} _{j\in\mathbb{N}}$ form an infinite sequence which goes to infinity, viz.
\[
0\le\mu_{0}<\mu_{1}\le\mu_{2}\le\ldots,\;\text{and }\lim_{j\to\infty}\mu_{j}=\infty.
\]
As a consequence, one can prove that $\left\langle \xi_{1},\phi_{j}\right\rangle $ and $\left\langle \xi_{2},\phi_{j}\right\rangle $
for $j=\overline{1,n}$ are i.i.d. standard Gaussian random variables. 

\subsection{QR-based methods and our novelty}\label{sec:1.3}
The QR method has a long remarkable history since the pioneering monograph \cite{MR0232549}. During the development of the QR method for inverse PDEs with deterministic noise, people usually focus on
\begin{itemize}
	\item the spectral methods that allow us to consider the solution in a mild presentation (cf., e.g., \cite{Tuan2017,Doan2017,Kaltenbacher2019}), but this typical method seems hard to handle the practical error control since one struggles with the Fourier accumulation for the nonlinear contexts;
	\item the Carleman-type estimate (cf., e.g., \cite{Klibanov2015,Klibanov2015a}) that includes the uniqueness result and convergence analysis, but currently, it only works for the short time observation with the linear equation.
\end{itemize}
Cf. \cite{Nguyen2019}, our approach removes this fence by relying on the so-called perturbing and stabilized operators, benefited from the essential bound of the nonlinear diffusion coefficient.
\begin{definition}[perturbing operator]\label{def:op1} The linear mapping $\mathbf{Q}_{\varepsilon}^{\beta}:[L^{2}(\Omega)]^{N}\to [L^{2}(\Omega)]^{N}$ is said to be a perturbing operator if there exist a function space $\mathbb{W}\subset [L^{2}(\Omega)]^{N}$ and an $\varepsilon$-independent constant $C_0>0$ such that 
	\[
	\left\Vert \mathbf{Q}_{\varepsilon}^{\beta}u\right\Vert _{\left[L^{2}\left(\Omega\right)\right]^{N}}\le\frac{C_{0}}{\gamma\left(T,\beta\right)}\left\Vert u\right\Vert _{\mathbb{W}}\quad\text{for any }u\in\mathbb{W}.
	\]
\end{definition}
\begin{definition}[stabilized operator]\label{def:op2} The linear mapping $\mathbf{P}_{\varepsilon}^{\beta}:[L^{2}(\Omega)]^{N}\to [L^{2}(\Omega)]^{N}$ is said to be a stabilized operator if there exists an $\varepsilon$-independent constant $C_1>0$ such that
	\[
	\left\Vert \mathbf{P}_{\varepsilon}^{\beta}u\right\Vert _{\left[L^{2}\left(\Omega\right)\right]^{N}}\le C_{1}\log\left(\gamma\left(T,\beta\right)\right)\left\Vert u\right\Vert _{\left[L^{2}\left(\Omega\right)\right]^{N}}\quad\text{for any }u\in\left[L^{2}\left(\Omega\right)\right]^{N}.
	\]
\end{definition}

In \cref{def:op1,def:op2}, $N$ denotes the number of species involved in the model and in this work, $N = 2$. Meanwhile, we denote by $\beta=\beta(\varepsilon)\in (0,1)$ the \emph{regularization parameter} satisfying $\lim_{\varepsilon\to 0^{+}}\beta(\varepsilon)=0$. The function $\gamma:[0,T]\times (0,1)\to \mathbb{R}$ indicates the decay behaviour of the perturbing operator, which stems from the source condition measuring the high smoothness of the true solution. In our analysis, we require: for any $\beta>0$, there holds
\begin{align}\label{abc}
\gamma\left(T,\beta\right)\ge1,\quad\lim_{\beta\to0^{+}}\gamma\left(t,\beta\right)=\infty\quad\text{for all }t\in\left(0,T\right].
\end{align}

In this study, we aim at developing this new QR approach and estimating its quality in a statistical setting. Cf. \cite{Alquier2011} the problem is well known to be exponentially ill-posed (compared to those mildly ill-posed introduced in the same reference); therefore, it is natural that we start from a very standard noise process. Yet, it is still a theoretical study that requires knowledge of the magnitude of the noise $\varepsilon$. Our theoretical analysis in this work and in all previous studies relies very much on the noise level as we want the (stably) approximate solution is close to the true one when $\varepsilon$ gets smaller. Without this analysis, there is nothing to ensure that the QR method performs well in solving such a highly challenging nonlinear inverse problem. Besides, this work would prepare a playground for the future evolution of our method. As we want to approach the real-world applications, it will concern more practical frameworks.

Although it is clear that in the framework of deterministic noise we can prove the error estimates for the QR method, it does not ensure that one can adapt those to the statistical inverse model in this work, including the references cited above. Besides, the existing literature on regularization of the system \eqref{system1} is very limited due to the inception stage. Starting off with the linear version of \eqref{system1} (where $\mathcal{D}_{i}$, $F$ and $G$ are independent of the solutions), we are aware of the presence of  \cite{Minh2018} where the authors exploited the trigonometric method in nonparametric regression to make use of the cut-off regularization in the statistical setting. This method is, however, very challenging in the nonlinear perspective as postulated in, e.g., \cite{Kirane2017}. It is worth noting in \cite{Kirane2017} that the authors essentially applied the QR-based method, where the Gevrey source condition is taken into account, to regularize a semilinear parabolic problem. Henceforth, our contribution herein can also be seen as an improvement of \cite{Kirane2017}, since the nonlocal spatial operator becomes rather challenging from the numerical standpoint in \cite{Kirane2017}.

From now on, some impedimenta that we will meet in proofs of our analysis should be revealed, except what we were very much concerned in \cite{Nguyen2019}. Due to the white noise processes we are taking into account, our observations, i.e. the measured terminal data, do not belong to $L^2(\Omega)$, but acts on $L^2(\Omega)$ by virtue of \eqref{eq:data1}--\eqref{eq:data2}. In principle, the random noise is large compared to the deterministic one. Thus, the adaptation of \cite{Nguyen2019} to this statistical scenario is not straightforward. Aside from the high smoothness of the true solution, we need the true terminal function to be very smooth to obtain rates of convergence in expectation.

\subsection{Outline of the paper}

The rest of the paper is organized as follows. In \cref{sec:2}, we introduce the notation and working assumptions for our analysis below. We also provide a ``computable'' example to validate the presence of the perturbation in \cref{def:op1}. In addition, we delineate a way to reconstruct the terminal data from the stochastic setting under consideration. Our main results are reported in \cref{sec:main}. Based on the proposed QR framework, we establish a regularization for the time-reversed system \eqref{system1}; cf. \cref{sec:3.1}. Convergence rates of the scheme are thoroughly explored in \cref{sec:error}; cf. \cref{thm:1,thm:2}. Proofs of these error estimates are detailed in \cref{sec:321,sec:proof2}, respectively. We close this main section by several discussions in \cref{sec:remarks}. Finally, two numerical tests are provided in \cref{sec:numerical} to verify our convergence analysis, and the conclusions follow in \cref{sec:conclusions}.

\section{Preliminaries}\label{sec:2}
In the sequel, wherever $\left\langle \cdot,\cdot\right\rangle $ and $\left\Vert \cdot\right\Vert $ are present, we mean the $L^2$ inner product and its corresponding norm. Meanwhile, the other standard Sobolev norms will be specified, if used. Also, we use $\mathcal{D}_{i}(u)(t)$ as $\mathcal{D}_{i}(\ell_{0}(u)(t))$  for ease of presentation. In this statistical inverse problem, we assume that the true solution exists uniquely and its regularity is assumed in \cref{thm:1,thm:2}. The backward uniqueness result for this problem will be studied in the near future. To this end, we use the following assumptions:

$\left(\text{A}_{1}\right)$ For $i=1,2$, the measurable functions
$\mathcal{D}_{i}>0$ is such that the mapping $\xi\mapsto\mathcal{D}_{i}\left(\xi\right)$
is continuous for $\xi\in\mathbb{R}$. Moreover, there exist $\varepsilon$-independent
$\underline{M},M_1,\overline{M}>0$ such that
\[
\underline{M}\le\mathcal{D}_{i}\left(\xi\right) \leq M_1 <\overline{M}\quad\text{for any }\xi\in\mathbb{R}.
\]

$\left(\text{A}_{2}\right)$ The source functions $F,G$ are measurable
and locally Lipschitz-continuous in the sense that
\begin{align*}
	& \left|F\left(x,t;u_{1};v_{1}\right)-F\left(x,t;u_{2};v_{2}\right)\right|+\left|G\left(x,t;u_{1};v_{1}\right)-G\left(x,t;u_{2};v_{2}\right)\right|\\
	& \le L\left(\ell\right)\left(\left|u_{1}-u_{2}\right|+\left|v_{1}-v_{2}\right|\right),
\end{align*}
for $\max\left\{ \left|u_{1}\right|,\left|u_{2}\right|,\left|v_{1}\right|,\left|v_{2}\right|\right\} \le\ell$
for some $\ell>0$.

$\left(\text{A}_{3}\right)$ The true final conditions $u_{f},v_{f}$ belong to $H^{2p}(\Omega)$ for $p>0$.

$\left(\text{A}_{1}'\right)$ For $i=1,2$, the diffusion $\mathcal{D}_{i}$ is globally Lipschitz-continuous in the sense that there exists $\varepsilon$-independent $\tilde{L}$ such that
\[
\left|\mathcal{D}_{i}\left(\ell_{0}\left(u_{1}\right)\left(t\right)\right)-\mathcal{D}_{i}\left(\ell_{0}\left(u_{2}\right)\left(t\right)\right)\right|\le\tilde{L}\left\Vert u_{1}\left(\cdot,t\right)-u_{2}\left(\cdot,t\right)\right\Vert \quad\text{for }u_{1},u_{2}\in L^{2}\left(\Omega\right).
\]

Since our observations $u_{f}^{\varepsilon}$ and $v_{f}^{\varepsilon}$ are not in general elements of $L^2(\Omega)$, we need the presence of $\left(\text{A}_{3}\right)$ to handle the following lemma (cf. \cite{Nane2017}):

	\begin{lemma}[Truncated Fourier reconstruction]\label{lem:2.1}
	Given $n\in \mathbb{N}^{*}$ the number of steps of discrete observations, we define 
	\[
	L^2(\Omega)\ni U_{f}^{\varepsilon,n}\left(x\right):=\sum_{j=1}^{n}\left\langle u_{f}^{\varepsilon},\phi_{j}\right\rangle \phi_{j}\left(x\right),\;
	L^2(\Omega)\ni V_{f}^{\varepsilon,n}\left(x\right):=\sum_{j=1}^{n}\left\langle v_{f}^{\varepsilon},\phi_{j}\right\rangle \phi_{j}\left(x\right). 
	\]
	Then with the aid of $\left(\text{A}_{3}\right)$, one has
	\[
	\mathbb{E}\left\Vert U_{f}^{\varepsilon,n}-u_{f}\right\Vert ^{2}\le\varepsilon^{2}n+\frac{\left\Vert u_{f}\right\Vert _{H^{2p}\left(\Omega\right)}^{2}}{\mu_{n}^{2p}},\quad\mathbb{E}\left\Vert V_{f}^{\varepsilon,n}-v_{f}\right\Vert ^{2}\le\varepsilon^{2}n+\frac{\left\Vert v_{f}\right\Vert _{H^{2p}\left(\Omega\right)}^{2}}{\mu_{n}^{2p}},
	\]
	where $\mu_{n}>0$ is the $n$th Dirichlet eigenvalue of the Laplacian operator.
\end{lemma}

The existence of the perturbing and stabilized operators (for \cref{def:op1,def:op2}) were commenced in \cite{Nguyen2019}, mimicking the stochastic gradient descent algorithm to obtain the function space $\mathbb{W}$ as a Gevrey\footnote{See again the Gevrey-like space defined in \cite[Section 5.2]{Nguyen2019}.} class of real-analytic functions. Below, we exemplify another one using the truncated Fourier method.

%

\begin{example}\label{example1}
	Upon the presence of the Dirichlet eigen-elements of the Laplacian operator, we choose
	\[
	\mathbf{Q}_{\varepsilon}^{\beta}u=\overline{M}\sum_{\mu_{j}>\frac{1}{\overline{M}T}\log\left(\gamma\left(T,\beta\right)\right)}\mu_{j}\left\langle u,\phi_{j}\right\rangle \phi_{j}\quad\text{for }u\in L^{2}\left(\Omega\right).
	\]
	Recall the Gevrey type class of functions of order $q>0$ and index $p>0$ defined by the spectrum of the Laplacian as follows:
	\[
	\mathbb{G}_{p,q}=\left\{ u\in L^{2}\left(\Omega\right):\sum_{j\in\mathbb{N}}\mu_{j}^{q}e^{2p\mu_{j}}\left|\left\langle u,\phi_{j}\right\rangle \right|^{2}<\infty\right\}.
	\]
	This is a Hilbert space equipped with the following inner product and norm:
	\begin{align*}
	& \left\langle u_{1},u_{2}\right\rangle _{\mathbb{G}_{p,q}}=\left\langle \left(-\Delta\right)^{q/2}e^{p\sqrt{-\Delta}}u_{1},\left(-\Delta\right)^{q/2}e^{p\sqrt{-\Delta}}u_{2}\right\rangle ,\\
	& \left\Vert u_{1}\right\Vert _{\mathbb{G}_{p,q}}=\sqrt{\sum_{j\in\mathbb{N}}\mu_{j}^{q}e^{2p\mu_{j}}\left|\left\langle u,\phi_{j}\right\rangle \right|^{2}}.
	\end{align*}
	With the help of the Parseval identity, it is then easy to see that $\bar{\mathbb{W}}=\mathbb{G}_{\overline{M}T,2}$ and $\bar{C}_0=\overline{M}$ in \cref{def:op1}. Now taking $\mathbf{P}_{\varepsilon}^{\beta}:=\overline{M}\Delta + \mathbf{Q}_{\varepsilon}^{\beta}$ to absorb high frequencies in the Laplacian, we get
	\[
	\mathbf{P}_{\varepsilon}^{\beta}u=-\overline{M}\sum_{\mu_{j}\le\frac{1}{\overline{M}T}\log\left(\gamma\left(T,\beta\right)\right)}\mu_{j}\left\langle u,\phi_{j}\right\rangle \phi_{j},
	\]
	and therefore, it holds that $C_1 = 1/T$ in \cref{def:op2}.
\end{example}

\begin{remark}\label{rem:1}
	The conventional cut-off function for the locally Lipschitz $F$ can be taken by
	\[
	F_{\ell^{\varepsilon}}\left(x,t;u,v\right):=\begin{cases}
	F\left(x,t;\ell^{\varepsilon};\ell^{\varepsilon}\right) & \text{if }\max\left\{ u,v\right\} >\ell^{\varepsilon},\\
	F\left(x,t;u;v\right) & \text{if }\max\left\{ u,v\right\} \in\left[-\ell^{\varepsilon},\ell^{\varepsilon}\right],\\
	F\left(x,t;-\ell^{\varepsilon};-\ell^{\varepsilon}\right) & \text{if }\max\left\{ u,v\right\} <-\ell^{\varepsilon},
	\end{cases}
	\]
	where $\ell^{\varepsilon}=\ell(\varepsilon)>0$ satisfying $\lim_{\varepsilon\to 0^{+}}\ell^{\varepsilon} = \infty$ is called the cut-off parameter. In the same manner, we have the cut-off function $G_{\ell^{\varepsilon}}(x,t;u;v)$ and thus one can prove that
	\begin{align*}
	&\left|F_{\ell^{\varepsilon}}\left(x,t;u_{1};v_{1}\right)-F_{\ell^{\varepsilon}}\left(x,t;u_{2};v_{2}\right)\right|+
	\left|G_{\ell^{\varepsilon}}\left(x,t;u_{1};v_{1}\right)-G_{\ell^{\varepsilon}}\left(x,t;u_{2};v_{2}\right)\right|\\
	&\le 2L\left(\ell^{\varepsilon}\right)\left(\left|u_{1}-u_{2}\right|+\left|v_{1}-v_{2}\right|\right),
	\end{align*}
	and one can observe that $\lim_{\varepsilon\to 0^{+}} L(\ell^{\varepsilon}) = \infty$. The proof can be detailed in \cite{Nguyen2019,Tuan2017}.
\end{remark}


\section{Main results}
\label{sec:main}

\subsection{Approximation by the QR method}\label{sec:3.1}
This part is devoted to establishing an approximate problem of \eqref{system1} by using the QR method. We begin by adding the perturbing operator to the original system as follows:
\[
\begin{cases}
u_{t}^{\varepsilon,n}-\mathcal{D}_{1}\left(u^{\varepsilon,n}\right)\Delta u^{\varepsilon,n}-\mathbf{Q}_{\varepsilon}^{\beta}u^{\varepsilon,n}=F\left(x,t;u^{\varepsilon,n};v^{\varepsilon,n}\right),\\
v_{t}^{\varepsilon,n}-\mathcal{D}_{2}\left(v^{\varepsilon,n}\right)\Delta v^{\varepsilon,n}-\mathbf{Q}_{\varepsilon}^{\beta}v^{\varepsilon,n}=G\left(x,t;u^{\varepsilon,n};v^{\varepsilon,n}\right) & \text{for }\left(x,t\right)\in Q_{T}.
\end{cases}
\]
Hereby, we define new diffusion-like coefficients $\overline{\mathcal{D}}_{i}$ for $i=1,2$ in such a way that $\overline{\mathcal{D}}_{i}:=\overline{M}-\mathcal{D}_{i}$. This way one can prove that $\overline{\mathcal{D}}_{i}\in (\overline{M}-M_1,\overline{M})$ for $M_1 < \overline{M}$ by $\left(\text{A}_{1}\right)$. Thus, we rely on the choice of the stabilized operator (i.e. $\mathbf{P}_{\varepsilon}^{\beta}:=\overline{M}\Delta + \mathbf{Q}_{\varepsilon}^{\beta}$) to arrive at
\begin{align}\label{eq:regularized}
\begin{cases}
u_{t}^{\varepsilon,n}+\overline{\mathcal{D}}_{1}\left(u^{\varepsilon,n}\right)\Delta u^{\varepsilon,n}=F_{\ell^{\varepsilon}}\left(x,t;u^{\varepsilon,n};v^{\varepsilon,n}\right) + \mathbf{P}_{\varepsilon}^{\beta}u^{\varepsilon,n},\\
v_{t}^{\varepsilon,n}+\overline{\mathcal{D}}_{2}\left(v^{\varepsilon,n}\right)\Delta v^{\varepsilon,n}=G_{\ell^{\varepsilon}}\left(x,t;u^{\varepsilon,n};v^{\varepsilon,n}\right)+\mathbf{P}_{\varepsilon}^{\beta}v^{\varepsilon,n} & \text{for }\left(x,t\right)\in Q_{T},
\end{cases}
\end{align}
where $F_{\ell^{\varepsilon}}$ and $G_{\ell^{\varepsilon}}$ are defined as in \cref{rem:1}. This is the coupled regularized system of \eqref{system1} that we wish to scrutinize in this work. Due to \cref{lem:2.1}, we associate \eqref{eq:regularized} with the Dirichlet boundary condition and the following terminal conditions:
\begin{align}\label{eq:regudata}
u^{\varepsilon,n}\left(x,T\right)=U_{f}^{\varepsilon,n}\left(x\right),\quad v^{\varepsilon,n}\left(x,T\right)=V_{f}^{\varepsilon,n}\left(x\right).
\end{align}
Hence, these equations form our regularized problem $\left(P^{\varepsilon,n}\right)$.

To study the weak solvability of the regularized problem $\left(P^{\varepsilon,n}\right)$, we use an exponential weight function $e^{\rho_{\varepsilon}(t-T)}$, where $\rho_{\varepsilon}>0$ is called as an $\varepsilon$-dependent auxiliary parameter, to consider the mappings $U^{\varepsilon,n}=e^{\rho_{\varepsilon}(t-T)}u^{\varepsilon,n}$ and $V^{\varepsilon,n}=e^{\rho_{\varepsilon}(t-T)}v^{\varepsilon,n}$. Thus, 
\eqref{eq:regularized} becomes
\begin{align}\label{eq:regunew}
\begin{cases}
U_{t}^{\varepsilon,n}+\overline{\mathcal{D}}_{1}\left(e^{\rho_{\varepsilon}\left(T-t\right)}U^{\varepsilon,n}\right)\Delta U^{\varepsilon,n}-\rho_{\varepsilon}U^{\varepsilon,n}\\
\qquad =e^{\rho_{\varepsilon}\left(t-T\right)}F_{\ell^{\varepsilon}}\left(e^{\rho_{\varepsilon}\left(T-t\right)}U^{\varepsilon,n};e^{\rho_{\varepsilon}\left(T-t\right)}V^{\varepsilon,n}\right)+\mathbf{P}_{\varepsilon}^{\beta}U^{\varepsilon,n},\\
V_{t}^{\varepsilon,n}+\overline{\mathcal{D}}_{2}\left(e^{\rho_{\varepsilon}\left(T-t\right)}V^{\varepsilon,n}\right)\Delta V^{\varepsilon,n}-\rho_{\varepsilon}V^{\varepsilon,n}\\
\qquad=e^{\rho_{\varepsilon}\left(t-T\right)}G_{\ell^{\varepsilon}}\left(e^{\rho_{\varepsilon}\left(T-t\right)}U^{\varepsilon,n};e^{\rho_{\varepsilon}\left(T-t\right)}V^{\varepsilon,n}\right)+\mathbf{P}_{\varepsilon}^{\beta}V^{\varepsilon,n},
\end{cases}
\end{align}
endowed with the Dirichlet boundary condition and with the same terminal data as \eqref{eq:regudata}. Henceforward, we define a weak formulation of this transformed system in the following type.

\begin{definition}\label{def:weak}
	For each $\varepsilon>0$, a pair of functions $\left(U^{\varepsilon,n},V^{\varepsilon,n}\right)$ is said to be
	a weak solution of \eqref{eq:regunew} if
	\[
	\left(U^{\varepsilon,n},V^{\varepsilon,n}\right)\in\left[L^{2}\left(0,T;H_{0}^{1}\left(\Omega\right)\right)\cap L^{\infty}\left(0,T;L^{2}\left(\Omega\right)\right)\right]^{2}
	\]
	and it holds that
	\begin{align}\label{eq:weak1}\frac{d}{dt}\left\langle  U^{\varepsilon,n},\psi_{1}\right\rangle  & -\overline{\mathcal{D}}_{1}\left(e^{\rho_{\varepsilon}\left(T-t\right)}U^{\varepsilon,n}\right)\int_{\Omega}\nabla U^{\varepsilon,n}\cdot\nabla\psi_{1}dx-\rho_{\varepsilon}\left\langle U^{\varepsilon,n},\psi_{1}\right\rangle \\
		\nonumber & =e^{\rho_{\varepsilon}\left(t-T\right)}\left\langle F_{\ell^{\varepsilon}}\left(e^{\rho_{\varepsilon}\left(T-t\right)}U^{\varepsilon,n};e^{\rho_{\varepsilon}\left(T-t\right)}V^{\varepsilon,n}\right),\psi_{1}\right\rangle +\left\langle \mathbf{P}_{\varepsilon}^{\beta}U^{\varepsilon,n},\psi_{1}\right\rangle ,\\
		\label{eq:weak2}\frac{d}{dt}\left\langle  V^{\varepsilon,n},\psi_{2}\right\rangle  & -\overline{\mathcal{D}}_{2}\left(e^{\rho_{\varepsilon}\left(T-t\right)}V^{\varepsilon,n}\right)\int_{\Omega}\nabla V^{\varepsilon,n}\cdot\nabla\psi_{2}dx-\rho_{\varepsilon}\left\langle V^{\varepsilon,n},\psi_{2}\right\rangle \\
		\nonumber & =e^{\rho_{\varepsilon}\left(t-T\right)}\left\langle G_{\ell^{\varepsilon}}\left(e^{\rho_{\varepsilon}\left(T-t\right)}U^{\varepsilon,n};e^{\rho_{\varepsilon}\left(T-t\right)}V^{\varepsilon,n}\right),\psi_{2}\right\rangle +\left\langle \mathbf{P}_{\varepsilon}^{\beta}V^{\varepsilon,n},\psi_{2}\right\rangle ,
	\end{align}
for all $\psi_1, \psi_2 \in H_0^1(\Omega)$.
\end{definition}

It is worth mentioning that to be successful with the Galerkin-type aid we choose in \cite{Nguyen2019}, the weight $\rho_{\varepsilon}$ must be large, controlled by the largeness of the magnitude stability of the regularized problem (cf. \cite{Nam2010}) and of the cut-off parameter $\ell^{\varepsilon}$. In the following, we only state the well-posedness result of the ``scaled" problem \eqref{eq:regunew}, while details of proof can be deduced as in \cite{Nguyen2019}.

\begin{theorem}
	Assume $\left(\text{A}_{1}\right)$--$\left(\text{A}_{3}\right)$ hold. For each $\varepsilon>0$, the regularized problem \eqref{eq:regunew} admits a pair of weak solutions $\left(U^{\varepsilon,n},V^{\varepsilon,n}\right)$ in the sense of \cref{def:weak}. Moreover, one has  $U^{\varepsilon,n},V^{\varepsilon,n}\in C([0,T];L^2(\Omega))$ and $U^{\varepsilon,n}_{t},V^{\varepsilon,n}_{t}\in L^2(0,T;(H^1(\Omega))')$.
\end{theorem}

In general, one can get the strong solution of \eqref{eq:regunew} in the sense that $U^{\varepsilon,n}_{t},V^{\varepsilon,n}_{t}$ and $\Delta U^{\varepsilon,n},\Delta V^{\varepsilon,n}$ belong to $L^2(0,T;L^2(\Omega))$ by increasing the regularity of the corresponding terminal conditions $U^{\varepsilon,n}_{f},V^{\varepsilon,n}_{f}$ (cf. \eqref{eq:regudata}) in $H^1(\Omega)$; see \cite[Theorem 3.3]{Almeida2016} for detailed techniques that we can adapt. However, this augment is not practical as in the context we are dealing with the noisy data. In other words, their gradients cannot be measured and even if it is possible, it eventually requires very much effort and expense. Accordingly, this explains why at present, we cannot fully adapt the error estimates for the finite element solution of \eqref{eq:regunew} obtained in, e.g., \cite{Duque2016,Chaudhary2018} in \cref{sec:numerical}. This open question will be explored in the near future.  

\subsection{Error analysis}\label{sec:error}
In the following, $n$ will be dependent of $\varepsilon$ (i.e. $n:=n(\varepsilon)\in\mathbb{N}$) due to the argument obtained in \cref{lem:2.1}, showing that $n$ cannot be arbitrarily large. Unlike the previous section where we design the approximate problem, here we denote by $u_{\beta}^{\varepsilon,n},v_{\beta}^{\varepsilon,n}$ the regularized solutions due to the involved regularization parameter $\beta$, recalled from \cref{sec:1.3}. Note here that the weight $\rho_{\beta}$ we use in proofs of the main results plays the same role as the weight $\rho_{\varepsilon}$ considered in \cref{sec:3.1}, i.e. it must be large as driven by the smallness of $\varepsilon$. Additionally, we below assume that $C_1 T \le 1$ and
\begin{align}\label{eq:gammaassume}
	\lim_{\varepsilon\to 0^{+}}\gamma^{C_{1}T}(T,\beta)\varepsilon\sqrt{n} = K_0\in (0,\infty),\;
	\lim_{\varepsilon\to 0^{+}}\gamma^{C_{1}T}(T,\beta)\mu_{n}^{-p} = K_1\in (0,\infty).
\end{align}

\begin{theorem}[error estimate for $0<t< T$]\label{thm:1} Assume $\left(\text{A}_{1}\right)$--$\left(\text{A}_{3}\right)$ and $\left(\text{A}_{1}'\right)$ hold. Consider the coupled Dirichlet system \eqref{system1} with terminal data \eqref{eq:noise} satisfying \eqref{eq:data1}--\eqref{eq:data2}. Suppose that it has a unique pair of solutions satisfying the source condition
	\begin{align*}
		u,v\in C([0,T];L^2(\Omega))\cap L^2(0,T;\bar{\mathbb{W}})\cap L^{\infty}(0,T;H^1_{0}(\Omega)\cap L^{\infty}(\Omega)),
	\end{align*}
where the function space $\bar{\mathbb{W}}$ is obtained from the choice of the perturbing operator $\mathbf{Q}_{\varepsilon}^{\beta}$ in \cref{def:op1}. By the resulting stabilized operator $\mathbf{P}_{\varepsilon}^{\beta}$ in \cref{def:op2}, we consider $\left(u^{\varepsilon,n}_{\beta},v^{\varepsilon,n}_{\beta}\right)$ as a pair of solutions of the approximate system \eqref{eq:regularized}--\eqref{eq:regudata}. Then for $0<\kappa\le 2C_1 t$ the following estimate holds:
\begin{align*}
	& \mathbb{E}\left(\left\Vert u_{\beta}^{\varepsilon,n}\left(\cdot,t\right)-u\left(\cdot,t\right)\right\Vert^2 +\left\Vert v_{\beta}^{\varepsilon,n}\left(\cdot,t\right)-v\left(\cdot,t\right)\right\Vert^2 \right)\\
	& +\frac{\overline{M}-M_1}{2}\mathbb{E}\int_{t}^{T}\left(\left\Vert \nabla\left(u_{\beta}^{\varepsilon,n}-u\right)\left(\cdot,s\right)\right\Vert^2 +\left\Vert \nabla\left(v_{\beta}^{\varepsilon,n}-v\right)\left(\cdot,s\right)\right\Vert^2 \right)ds\\
	& \le\left[K_{0}^{2}+\left(K_{1}^{2}+\bar{C}_{0}^{2}\right)\tilde{C}_{p}\right]\gamma^{-2C_{1}t}\left(T,\beta\right)\log^{\kappa}\left(\gamma(T,\beta)\right)e^{2\left(T-t\right)C_{3}},
\end{align*}
where $C_3>0$, $\tilde{C}_p>0$ are not dependent of $\varepsilon$. 
\end{theorem}

There is no doubt that it is hard to attain the convergence of regularization at $t=0$. Herewith, by requiring more information of the source condition we find $t^{\varepsilon}>0$ in such a way that $u_{\beta}^{\varepsilon}(\cdot,t=t^{\varepsilon})$ will be a good approximation of $u(\cdot,t=0)$, involving the expectation operator.

\begin{theorem}[error estimate for $t=0$]\label{thm:2}
	Under the assumptions of \cref{thm:1}, we further assume that
	$u,v\in C^1(0,T;L^2(\Omega))$. Then there exists $t^{\varepsilon}\in (0,T)$ approaching 0 such that
	\begin{align*}
		& \mathbb{E}\left\Vert u_{\beta}^{\varepsilon,n}\left(\cdot,t^{\varepsilon}\right)-u\left(\cdot,0\right)\right\Vert ^{2}+\mathbb{E}\left\Vert v_{\beta}^{\varepsilon,n}\left(\cdot,t^{\varepsilon}\right)-v\left(\cdot,0\right)\right\Vert ^{2}\\
		& \le\left[K_{0}^{2}+\left(K_{1}^{2}+\bar{C}_{0}^{2}\right)\tilde{C}_{p}\right]C_{1}^{-1}\log^{\kappa^{\varepsilon}-1}\left(\gamma\left(T,\beta\right)\right)e^{2TC_{3}} +C_{1}^{-1}\log^{-1}\left(\gamma\left(T,\beta\right)\right)C_{4},
	\end{align*}
where $C_4 > 0$ is independent of $\varepsilon$ and $0<\kappa^{\varepsilon}<\min\left\{2C_1 t^{\varepsilon},1\right\}$.
\end{theorem}

\subsubsection{Proof of \cref{thm:1}}\label{sec:321}
Put $\mathcal{X}_{\beta}^{\varepsilon,n}\left(x,t\right)=e^{\rho_{\beta}\left(t-T\right)}\left[u_{\beta}^{\varepsilon,n}\left(x,t\right)-u\left(x,t\right)\right]$,
$\mathcal{Y}_{\beta}^{\varepsilon,n}\left(x,t\right)=e^{\rho_{\beta}\left(t-T\right)}\left[v_{\beta}^{\varepsilon,n}\left(x,t\right)-v\left(x,t\right)\right]$ for $\rho_{\beta}>0$. Then the ``scaled" difference equation for $u$ is given by
\begin{align*}
	& \frac{\partial\mathcal{X}_{\beta}^{\varepsilon,n}}{\partial t}  +\overline{\mathcal{D}}_{1}\left(u_{\beta}^{\varepsilon,n}\right)\Delta\mathcal{X}_{\beta}^{\varepsilon,n}-\rho_{\beta}\mathcal{X}_{\beta}^{\varepsilon,n} =\mathbf{P}_{\varepsilon}^{\beta}\mathcal{X}_{\beta}^{\varepsilon,n}+e^{\rho_{\beta}\left(t-T\right)}\mathbf{Q}_{\varepsilon}^{\beta}u\\
	& -e^{\rho_{\beta}\left(t-T\right)}\left[\overline{\mathcal{D}}_{1}\left(u_{\beta}^{\varepsilon,n}\right)-\overline{\mathcal{D}}_{1}\left(u\right)\right]\Delta u +e^{\rho_{\beta}\left(t-T\right)}\left[F_{\ell^{\varepsilon}}\left(u_{\beta}^{\varepsilon};v_{\beta}^{\varepsilon}\right)-F\left(u;v\right)\right].
\end{align*}
In parallel, we also derive that the difference equation is endowed with the zero Dirichlet boundary condition $\mathcal{X}_{\beta}^{\varepsilon,n}=0$ on the boundary $\partial\Omega$ and with the terminal condition $\mathcal{X}_{\beta}^{\varepsilon,n}\left(x,T\right)=U_{f}^{\varepsilon,n}\left(x\right)-u_{f}\left(x\right)$
for $x\in\Omega$. Hereby, we multiply this equation by $\mathcal{X}_{\beta}^{\varepsilon,n}$ and then integrate the resulting equation over the domain of interest $\Omega$ to arrive at
\begin{align*}
	& \frac{1}{2}\frac{d}{dt}\left\Vert \mathcal{X}_{\beta}^{\varepsilon,n}\right\Vert ^{2}-\overline{\mathcal{D}}_{1}\left(u_{\beta}^{\varepsilon,n}\right)\left\Vert \nabla\mathcal{X}_{\beta}^{\varepsilon,n}\right\Vert ^{2}-\rho_{\beta}\left\Vert \mathcal{X}_{\beta}^{\varepsilon,n}\right\Vert ^{2} =\underbrace{\left\langle \mathbf{P}_{\varepsilon}^{\beta}\mathcal{X}_{\beta}^{\varepsilon,n},\mathcal{X}_{\beta}^{\varepsilon,n}\right\rangle }_{:=I_{1}}\\
	& +\underbrace{e^{\rho_{\beta}\left(t-T\right)}\left\langle \mathbf{Q}_{\varepsilon}^{\beta}u,\mathcal{X}_{\beta}^{\varepsilon,n}\right\rangle }_{:=I_{2}}+\underbrace{e^{\rho_{\beta}\left(t-T\right)}\left\langle F_{\ell^{\varepsilon}}\left(u_{\beta}^{\varepsilon,n};v_{\beta}^{\varepsilon,n}\right)-F\left(u;v\right),\mathcal{X}_{\beta}^{\varepsilon,n}\right\rangle }_{:=I_{3}}\\
	& +\underbrace{e^{\rho_{\beta}\left(t-T\right)}\left[\overline{\mathcal{D}}_{1}\left(u_{\beta}^{\varepsilon,n}\right)-\overline{\mathcal{D}}_{1}\left(u\right)\right]\left\langle \nabla u,\nabla\mathcal{X}_{\beta}^{\varepsilon,n}\right\rangle }_{:=I_{4}}.
\end{align*}

Cf. \cref{def:op2}, we bound $I_1$ from below by
\[
I_{1}\ge-C_1\log\left(\gamma\left(T,\beta\right)\right)\left\Vert \mathcal{X}_{\beta}^{\varepsilon,n}\right\Vert ^{2},
\]
aided by the Cauchy--Schwarz inequality. Meanwhile, we use the Young inequality in combination with \cref{def:op1} to get
\[
I_{2}\ge-\gamma^{-2}\left(T,\beta\right)e^{2\rho_{\beta}\left(t-T\right)}\bar{C}_{0}^{2}\left\Vert u\right\Vert _{\bar{\mathbb{W}}}^{2}-\frac{1}{4}\left\Vert \mathcal{X}_{\beta}^{\varepsilon,n}\right\Vert ^{2}.
\]

Now, by decreasing $\varepsilon$ we can choose the cut-off parameter $\ell^{\varepsilon}$ large in such a way that $\ell^{\varepsilon}\ge\max\left\{ \left\Vert u\right\Vert _{L^{\infty}\left(0,T;L^{2}\left(\Omega\right)\right)};\left\Vert v\right\Vert _{L^{\infty}\left(0,T;L^{2}\left(\Omega\right)\right)}\right\} $. Thus, it holds that $F_{\ell^{\varepsilon}}\left(u;v\right)=F\left(u;v\right)$, according to \cref{rem:1}. This way we bound $I_3$ from below by
\[
I_{3}\ge-4L\left(\ell^{\varepsilon}\right)\left(\left\Vert \mathcal{X}_{\beta}^{\varepsilon,n}\right\Vert ^{2}+\left\Vert \mathcal{Y}_{\beta}^{\varepsilon,n}\right\Vert ^{2}\right).
\]

Using $\left(\text{A}_{1}'\right)$, the estimate of $I_4$ is thus given by
\[
I_{4}\ge-\tilde{L}\left\Vert \mathcal{X}_{\beta}^{\varepsilon,n}\right\Vert \left\Vert \nabla u\right\Vert \left\Vert \nabla\mathcal{X}_{\beta}^{\varepsilon,n}\right\Vert \ge-\frac{\overline{M}-M_1}{4}\left\Vert \nabla\mathcal{X}_{\beta}^{\varepsilon,n}\right\Vert ^{2}-\frac{\tilde{L}^{2}\left\Vert \nabla u\right\Vert ^{2}}{\overline{M}-M_1}\left\Vert \mathcal{X}_{\beta}^{\varepsilon,n}\right\Vert ^{2}.
\]

Grouping the estimates of $I_i$ for $i=\overline{1,4}$, we obtain the fact that
\begin{align}\label{eq:3.7}
	& \left\Vert \mathcal{X}_{\beta}^{\varepsilon,n}\left(\cdot,T\right)\right\Vert ^{2}-\left\Vert \mathcal{X}_{\beta}^{\varepsilon,n}\left(\cdot,t\right)\right\Vert ^{2} + \rho_{\beta}^{-1}\gamma^{-2}\left(T,\beta\right)\left(1-e^{2\rho_{\beta}\left(t-T\right)}\right) \bar{C}_{0}^{2}\left\Vert u\right\Vert _{C(0,T;\bar{\mathbb{W}})}^{2}\\\nonumber
	& \ge\frac{\overline{M}-M_1}{2}\int_{t}^{T}\left\Vert \nabla\mathcal{X}_{\beta}^{\varepsilon,n}\left(\cdot,s\right)\right\Vert ^{2}ds-8L\left(\ell^{\varepsilon}\right)\int_{t}^{T}\left(\left\Vert \mathcal{X}_{\beta}^{\varepsilon,n}\left(\cdot,s\right)\right\Vert ^{2}+\left\Vert \mathcal{Y}_{\beta}^{\varepsilon,n}\left(\cdot,s\right)\right\Vert ^{2}\right)ds\\\nonumber
	& +\left(2\rho_{\beta}-2C_{1}\log\left(\gamma\left(T,\beta\right)\right)-\frac{1}{2}-\frac{2\tilde{L}^{2}\left\Vert \nabla u\right\Vert_{L^{\infty}(0,T;L^2(\Omega))} ^{2}}{\overline{M}-M_1}\right)\int_{t}^{T}\left\Vert \mathcal{X}_{\beta}^{\varepsilon,n}\left(\cdot,s\right)\right\Vert ^{2}ds.
\end{align}

By doing the same token with the difference equation for $v$, we can find that
\begin{align}\label{eq:3.8}
	& \left\Vert \mathcal{Y}_{\beta}^{\varepsilon,n}\left(\cdot,T\right)\right\Vert ^{2}-\left\Vert \mathcal{Y}_{\beta}^{\varepsilon,n}\left(\cdot,t\right)\right\Vert ^{2}+ \rho_{\beta}^{-1}\gamma^{-2}\left(T,\beta\right)\left(1-e^{2\rho_{\beta}\left(t-T\right)}\right)\bar{C}_{0}^{2}\left\Vert v\right\Vert _{C(0,T;\bar{\mathbb{W}})}^{2}\\\nonumber
	& \ge\frac{\overline{M}-M_1}{2}\int_{t}^{T}\left\Vert \nabla\mathcal{Y}_{\beta}^{\varepsilon,n}\left(\cdot,s\right)\right\Vert ^{2}ds-8L\left(\ell^{\varepsilon}\right)\int_{t}^{T}\left(\left\Vert \mathcal{X}_{\beta}^{\varepsilon,n}\left(\cdot,s\right)\right\Vert ^{2} + \left\Vert \mathcal{Y}_{\beta}^{\varepsilon,n}\left(\cdot,s\right)\right\Vert ^{2}\right)ds\\\nonumber
	& +\left(2\rho_{\beta}-2C_{1}\log\left(\gamma\left(T,\beta\right)\right)-\frac{1}{2}-\frac{2\tilde{L}^{2}\left\Vert \nabla v\right\Vert_{L^{\infty}(0,T;L^2(\Omega))} ^{2}}{\overline{M}-M_1}\right)\int_{t}^{T}\left\Vert \mathcal{Y}_{\beta}^{\varepsilon,n}\left(\cdot,s\right)\right\Vert ^{2}ds.
\end{align}

Combining \eqref{eq:3.7}, \eqref{eq:3.8}, taking the expectation and recalling  \cref{lem:2.1}, we choose
\[
\rho_{\beta}= C_{1}\log\left(\gamma\left(T,\beta\right)\right) + \frac{1}{4} + 4L(\ell^{\varepsilon}) + \frac{\tilde{L}^{2}\left(\left\Vert \nabla u\right\Vert_{L^{\infty}(0,T;L^2(\Omega))} ^{2}+\left\Vert \nabla v\right\Vert_{L^{\infty}(0,T;L^2(\Omega))} ^{2}\right)}{\overline{M}-M_1}
\]
and thus, it yields
\begin{align*}
& \mathbb{E}\left(\left\Vert \mathcal{X}_{\beta}^{\varepsilon,n}\left(\cdot,t\right)\right\Vert ^{2}+\left\Vert \mathcal{Y}_{\beta}^{\varepsilon,n}\left(\cdot,t\right)\right\Vert ^{2}\right)+\frac{\overline{M}-M_1}{2}\mathbb{E}\int_{t}^{T}\left(\left\Vert \nabla\mathcal{X}_{\beta}^{\varepsilon,n}\left(\cdot,s\right)\right\Vert ^{2}+\left\Vert \nabla\mathcal{Y}_{\beta}^{\varepsilon,n}\left(\cdot,s\right)\right\Vert ^{2}\right)ds\\
& \le\varepsilon^{2}n+\left(\mu_{n}^{-2p}+\rho_{\beta}^{-1}\gamma^{-2}\left(T,\beta\right)\bar{C}_{0}^{2}\right)\tilde{C}_{p}.
\end{align*}
with $\tilde{C}_p>0$ given by
\begin{align*}
\tilde{C}_{p} & =\left\Vert u_{f}\right\Vert _{H^{2p}\left(\Omega\right)}^{2}+\left\Vert v_{f}\right\Vert _{H^{2p}\left(\Omega\right)}^{2}+\left\Vert u\right\Vert _{C\left(0,T;\bar{\mathbb{W}}\right)}^{2}\\
& +\left\Vert v\right\Vert _{C\left(0,T;\bar{\mathbb{W}}\right)}^{2}+\left\Vert \nabla u\right\Vert _{L^{\infty}\left(0,T;L^{2}\left(\Omega\right)\right)}^{2}+\left\Vert \nabla v\right\Vert _{L^{\infty}\left(0,T;L^{2}\left(\Omega\right)\right)}^{2}.
\end{align*}

By the back-substitutions of $\mathcal{X}_{\beta}^{\varepsilon,n}$ and $\mathcal{Y}_{\beta}^{\varepsilon,n}$, one has
\begin{align}\label{eq:3.9}
	& \mathbb{E}\left(\left\Vert u_{\beta}^{\varepsilon,n}\left(\cdot,t\right)-u\left(\cdot,t\right)\right\Vert ^{2}+\left\Vert v_{\beta}^{\varepsilon,n}\left(\cdot,t\right)-v\left(\cdot,t\right)\right\Vert ^{2}\right)\\\nonumber
	& +\frac{\overline{M}-M_1}{2}\mathbb{E}\int_{t}^{T}\left(\left\Vert \nabla\left(u_{\beta}^{\varepsilon,n}-u\right)\left(\cdot,s\right)\right\Vert ^{2}+\left\Vert \nabla\left(v_{\beta}^{\varepsilon,n}-v\right)\left(\cdot,s\right)\right\Vert ^{2}\right)ds\\\nonumber
	& \le\left[\varepsilon^{2}n+\left(\mu_{n}^{-2p}+\rho_{\beta}^{-1}\gamma^{-2}\left(T,\beta\right)\bar{C}_{0}^{2}\right)\tilde{C}_{p}\right]\gamma^{2C_{1}\left(T-t\right)}\left(T,\beta\right)e^{2\left(T-t\right)C_{2}\left(\ell^{\varepsilon}\right)},
\end{align}
where we have denoted by $C_{2}\left(\ell^{\varepsilon}\right)=\frac{1}{4}+4L(\ell^{\varepsilon})+\frac{\tilde{L}^{2}\tilde{C}_{p}}{\overline{M}-M_1}$ 

Observe in \eqref{eq:3.9} that if we take $n$ and $\gamma$ satisfying \eqref{eq:gammaassume} and choose $C_1$ such that $C_1 T \le 1$, we can go on the estimate \eqref{eq:3.9} as follows:
\begin{align}\label{eq:3.10}
	& \mathbb{E}\left(\left\Vert u_{\beta}^{\varepsilon,n}\left(\cdot,t\right)-u\left(\cdot,t\right)\right\Vert ^{2}+\left\Vert v_{\beta}^{\varepsilon,n}\left(\cdot,t\right)-v\left(\cdot,t\right)\right\Vert ^{2}\right)\\\nonumber
	& +\frac{\overline{M}-M_1}{2}\mathbb{E}\int_{t}^{T}\left(\left\Vert \nabla\left(u_{\beta}^{\varepsilon,n}-u\right)\left(\cdot,s\right)\right\Vert ^{2}+\left\Vert \nabla\left(v_{\beta}^{\varepsilon,n}-v\right)\left(\cdot,s\right)\right\Vert ^{2}\right)ds\\\nonumber
	& \le\left[K_{0}^{2}+\left(K_{1}^{2}+\bar{C}_{0}^{2}\right)\tilde{C}_{p}\right]\gamma^{-2C_{1}t}\left(T,\beta\right)e^{2\left(T-t\right)C_{2}\left(\ell^{\varepsilon}\right)}.
\end{align}

It now remains to deal with a fine control of $C_{2}(\ell^{\varepsilon})$ in \eqref{eq:3.10} that involves the blow-up profile of $L(\ell^{\varepsilon})$ (see again in \cref{rem:1}) as $\varepsilon\to 0^{+}$. In fact, due to the conventional logarithmic rate of convergence in regularization of parabolic problems, our attempt is to bound this $\varepsilon$-dependent constant in such a way that its growth does not destroy the aimed speed. To do that, we basically follow the strategy in \cite[Section 5]{Nguyen2019}, which enables us to choose
\begin{align}\label{eq:Lell}
L(\ell^{\varepsilon})\le \frac{1}{8(T-t)}\log\left(\log^{\kappa} (\gamma(T,\beta))\right)\quad \text{for }t\in (0,T),
\end{align}
where $\kappa>0$ is an $\varepsilon$-independent constant. Thus, it gives
$
e^{8(T-t)L(\ell^{\varepsilon})}\le \log^{\kappa}(\gamma(T,\beta))
$
and consequently, it holds that
\begin{align*}
	& \mathbb{E}\left(\left\Vert u_{\beta}^{\varepsilon,n}\left(\cdot,t\right)-u\left(\cdot,t\right)\right\Vert ^{2}+\left\Vert v_{\beta}^{\varepsilon,n}\left(\cdot,t\right)-v\left(\cdot,t\right)\right\Vert ^{2}\right)\\
	& +\frac{\overline{M}-M_1}{2}\mathbb{E}\int_{t}^{T}\left(\left\Vert \nabla\left(u_{\beta}^{\varepsilon,n}-u\right)\left(\cdot,s\right)\right\Vert ^{2}+\left\Vert \nabla\left(v_{\beta}^{\varepsilon,n}-v\right)\left(\cdot,s\right)\right\Vert ^{2}\right)ds\\
	& \le\left[K_{0}^{2}+\left(K_{1}^{2}+\bar{C}_{0}^{2}\right)\tilde{C}_{p}\right]\gamma^{-2C_{1}t}\left(T,\beta\right)\log^{\kappa}\left(\gamma(T,\beta)\right)e^{2\left(T-t\right)C_{3}}.
\end{align*}
where $C_3=\frac{1}{4}+\frac{\tilde{L}^{2}\tilde{C}_{p}}{\overline{M}-M_1}>0$ is now independent of $\varepsilon$.

Hence, by choosing $\kappa:=\kappa(t) \le 2C_1 t$ we complete the proof of the theorem.

\subsubsection{Proof of \cref{thm:2}}\label{sec:proof2}
Using the triangle inequality and aided by \cref{thm:1}, one has
\begin{align*}
	& \mathbb{E}\left\Vert u_{\beta}^{\varepsilon,n}\left(\cdot,t^{\varepsilon}\right)-u\left(\cdot,0\right)\right\Vert ^{2}+\mathbb{E}\left\Vert v_{\beta}^{\varepsilon,n}\left(\cdot,t^{\varepsilon}\right)-v\left(\cdot,0\right)\right\Vert ^{2}\\
	& \le2\mathbb{E}\left\Vert u_{\beta}^{\varepsilon}\left(\cdot,t^{\varepsilon}\right)-u\left(\cdot,t^{\varepsilon}\right)\right\Vert ^{2}+2\mathbb{E}\left\Vert v_{\beta}^{\varepsilon}\left(\cdot,t^{\varepsilon}\right)-v\left(\cdot,t^{\varepsilon}\right)\right\Vert ^{2}\\
	& +2\mathbb{E}\left\Vert u\left(\cdot,t^{\varepsilon}\right)-u\left(\cdot,0\right)\right\Vert ^{2}+2\mathbb{E}\left\Vert v\left(\cdot,t^{\varepsilon}\right)-v\left(\cdot,0\right)\right\Vert ^{2}\\
	& \le\left[K_{0}^{2}+\left(K_{1}^{2}+\bar{C}_{0}^{2}\right)\tilde{C}_{p}\right]\gamma^{-2C_{1}t^{\varepsilon}}\left(T,\beta\right)\log^{\kappa^{\varepsilon}}\left(\gamma\left(T,\beta\right)\right)e^{2\left(T-t^{\varepsilon}\right)C_{3}}\\
	& +(t^{\varepsilon})^{2}\left(\left\Vert u_{t}\right\Vert _{C\left(0,T;L^{2}\left(\Omega\right)\right)}^{2}+\left\Vert v_{t}\right\Vert _{C\left(0,T;L^{2}\left(\Omega\right)\right)}^{2}\right).
\end{align*}

Here, $\kappa^{\varepsilon}:=\kappa(t^{\varepsilon})\in(0,2C_1 t^{\varepsilon}]$. Observe that since with respect to $t^{\varepsilon}$, the term $\gamma^{-C_1 t^{\varepsilon}}$ is decreasing, while $t^{\varepsilon}$ itself increases linearly. Thus, for each $\varepsilon > 0$ we can find a unique $t^{\varepsilon}\in (0,T)$ such that $t^{\varepsilon} = \gamma^{-C_1 t^{\varepsilon}}(T,\beta)$
or equivalently,
$
\frac{\log(t^{\varepsilon})}{t^{\varepsilon}}=-C_1 \log(\gamma(T,\beta))
$. Thanks to the inequality $\log(a)>-a^{-1}$ for all $a>0$, we find that
\[
t^{\varepsilon}<\sqrt{\frac{1}{C_1 \log(\gamma(T,\beta))}}.
\]
By the growth of $\gamma(T,\beta)$ as $\varepsilon\to 0^{+}$, one knows that $\lim_{\varepsilon\to 0^{+}} t^{\varepsilon} = 0$. Hence, we complete the proof of the theorem.

\subsubsection{Essential remarks}\label{sec:remarks}
Upon our main results reported in \cref{thm:1,thm:2}, we single out some important remarks where we believe that it is essential to highlight our novel QR method.
\begin{itemize}
	\item By enjoying the exponential weight function in proofs of these theorems, it is clear that the speeds of convergence are exponentially controlled with respect to time and this control becomes larger when our reconstruction runs close to $t=0$.
	\item When $F$ and $G$ are globally Lipschitz-continuous with respect to the arguments $u$ and $v$, we obtain the H\"older-type rate for $t\in (0,T)$, which means $\mathcal{O}\left(\gamma^{-2C_1t}(T,\beta)\right)$, whilst it is of the logarithmic order $\mathcal{O}(\log^{-1}(\gamma(T,\beta)))$ at $t=0$. These are what we have obtained in \cite[Theorem 4.5]{Nguyen2019}. The presence of the term $\log^{\kappa}(\gamma(T,\beta))$ in \cref{thm:1} shows that the locally Lipschitz-continuous nonlinearities usually encountered in real-world applications require much effort to approximate. In the near future, we will figure out a better way to design a new bound for the cut-off constant $L(\ell^{\varepsilon})$  (cf. \eqref{eq:Lell}).
	\item Cf. \cite{Nguyen2019}, if the model under consideration involves the gradient terms accounting for, e.g., haptotaxis in a population of cells (as presented in a model for haptotaxis and chemotaxis effects on cells' motion; cf. \cite{Sherratt1993}), our proposed method is still valid. A slightly lower rate of convergence is expected in controlling the gradient terms, but theoretically, it still converges logarithmically.
	\item If the condition \eqref{eq:gammaassume} is changed to  \begin{align}\label{eq:gammaassume1}
	\lim_{\varepsilon\to 0^{+}}\gamma^{C_{1}T}(T,\beta)\sqrt{\varepsilon n} = \widetilde{K}_0,\;
	\lim_{\varepsilon\to 0^{+}}\gamma^{C_{1}T}(T,\beta)\varepsilon^{-1/2}\mu_{n}^{-p} = \widetilde{K}_1,
	\end{align}
	we eventually obtain the convergence rate $\mathcal{O}\left(\varepsilon^{1/2} + \log^{-1}(\gamma(T,\beta))\right)$ in the space $C(t,T;L^2(\Omega))\cap L^2(t,T;H^1_0(\Omega))$ for any $0\le t<T$. This is easily deduced from the estimate \eqref{eq:3.9}. Therefore, the presence of \cref{thm:2} as well as the additional regularity assumption that $u,v\in C^1(0,T;L^2(\Omega))$ can be neglected. This is also a new finding compared to the convergence result in our original work \cite{Nguyen2019}.
	\item The special difference between the deterministic noise in \cite{Nguyen2019,Nam2010,Klibanov2015a,Tuan2017,Long1994} and the stochastic noise we are taking into account is the following. It is not fully about the involved expectation operator, but about the findings of $K_0$ and $K_1$ in \eqref{eq:gammaassume}. It shows that our QR method in this context really needs a very careful adaptation because not only $n$ has to be chosen properly, but also the eigenvalues driven by such $n$ need a fine control with respect to $\varepsilon$ and $\gamma$. Consider $\Omega$, for example, as an open parallelepiped $(0,a_1)\times\ldots(0,a_d)\subset \mathbb{R}^{d}$ with $a_{i}>0$, $i\in\left\{1,\ldots,d\right\}$. As one of advantages of this QR method, we only need to solve the simple Dirichlet eigenvalue problem regardless of the complex structure involved in the diffusion term, whenever it is essentially bounded. Henceforth, the Dirichlet eigen-elements are given by
	\[
	\phi_{l}=\prod_{j=1}^{d}\sqrt{\frac{2}{a_{j}}}\sin\left(\frac{\pi l_{j}}{a_{j}}x_{j}\right),\quad\mu_{l}=\sum_{j=1}^{d}\left(\frac{\pi l_{j}}{a_{j}}\right)^{2}\quad\text{for }l_{j}\in\mathbb{N},j\in\left\{ 1,\ldots,d\right\} .
	\]
	In this circumstance, we mean $n$ in \eqref{eq:gammaassume} as $\left|n\right|$ where $n:=\left(n_1,\ldots,n_{d}\right)$. Then, we can choose $n_j\propto\varepsilon^{-2\theta}$ for $j=1,\ldots,d$ and  $\theta\in (0,1)$. Notice that $\mu_{n}\propto\left|n\right|^2$, then $\gamma^{C_1 T}(T,\beta)\propto \varepsilon^{\max\left\{\theta-1,-4\theta p\right\}}$ is needed. Essentially, we have drawn a possibility to choose $n$ and $\gamma$ in the context of \eqref{eq:gammaassume}. Additionally, one can take $\beta = \varepsilon$ and thus $\gamma(t,\beta)=\varepsilon^{\frac{t}{C_1 T^2}\max\left\{\theta - 1, -4\theta p\right\}}$. Hence, the following error bounds hold
	\begin{align*}
		& \mathbb{E}\left\Vert u_{\beta}^{\varepsilon,n}\left(\cdot,t\right)-u\left(\cdot,t\right)\right\Vert ^{2}  +\mathbb{E}\left\Vert v_{\beta}^{\varepsilon,n}\left(\cdot,t\right)-v\left(\cdot,t\right)\right\Vert ^{2}\\
		& \le C\varepsilon^{\frac{2t}{T}\min\left\{ 1-\theta,4\theta p\right\} }\log^{\kappa(t)}\left(\varepsilon^{\frac{1}{C_1 T}\max\left\{\theta - 1,-4\theta p\right\}}\right)\quad\text{for }t\in\left(0,T\right),
	\end{align*}
	\begin{align*}
		&\mathbb{E}\left\Vert u_{\beta}^{\varepsilon,n}\left(\cdot,t^{\varepsilon}\right)-u\left(\cdot,0\right)\right\Vert ^{2}  +\mathbb{E}\left\Vert v_{\beta}^{\varepsilon,n}\left(\cdot,t^{\varepsilon}\right)-v\left(\cdot,0\right)\right\Vert ^{2}\\
		& \le C\log^{\kappa^{\varepsilon}-1}\left(\varepsilon^{\frac{1}{C_1 T}\max\left\{\theta - 1,-4\theta p\right\}}\right),
	\end{align*}
where $\kappa(t)\in(0,2C_1 t]$ and  $\kappa^{\varepsilon}\in\left(0,\min\left\{2C_1 t^{\varepsilon},1\right\}\right)$.
	\item In case of deterministic noises, the current analysis only needs the true terminal data in $L^2(\Omega)$. Moreover, the assumptions \eqref{eq:gammaassume} will be as we required in \cite[Theorem 4.5]{Nguyen2019}.
\end{itemize}


\section{Numerical tests}\label{sec:numerical}
Throughout the examples, we choose $\Omega = (a,b) = (0,\pi)$ as a sample domain of tissue in one-dimensional. Accordingly, the Dirichlet eigen-elements are elementary and have the following form:
\[
\phi_{j}(x) = \sqrt{\frac{2}{b-a}}\sin(jx),\quad \mu_{j} = j^2.
\]

Given $\varepsilon >0$, as discussed in \cref{sec:remarks} we choose that $n=\left\lfloor \varepsilon^{-2\theta}\right\rfloor$ for $\theta\in(0,1)$, $\beta(\varepsilon) = \varepsilon$ and $\gamma(t,\varepsilon) = \varepsilon^{\frac{t}{C_1 T^2}\max\left\{\theta - 1,-4\theta p \right\}}$ for all $t\in[0,T]$. Then \eqref{eq:gammaassume} is satisfied with $K_0,K_1\le 1$. Employing \cref{example1}, we define the stabilized operator $\mathbf{P}_{\varepsilon}^{\beta}$ as follows:
\[
\mathbf{P}_{\varepsilon}^{\beta}u=-\overline{M}\sum_{j\in\mathbb{B}_{\varepsilon}}\mu_{j}\left\langle u,\phi_{j}\right\rangle \phi_{j},
\]
where $\mathbb{B}_{\varepsilon}:= \left\{j\in\mathbb{N}:\mu_j \le \log\left(\varepsilon^{\frac{1}{C_1 T}\max\left\{\theta - 1,-4\theta p\right\}}\right) \right\}$ is the admissible set of our Fourier frequencies, and $C_1 = \overline{M}$.

We generate the discrete observations $u_{f}^{\varepsilon,j},v_{f}^{\varepsilon,j}$ for $1\le j \le n$ by the assumptions \eqref{eq:data1}--\eqref{eq:data2} where $\left\langle \xi_{i},\phi_j \right\rangle $ for $i=1,2$ are random variables with normal distribution.  In our illustration, we employ the manufactured solution to verify the validity of the regularization scheme. When doing so, we know the explicit expression of the true $u_{f}$ and $v_{f}$. Cf. \cite[Section 4]{Khoa2017}, the adaptive Filon-type method is rather reliable to compute the inner products $\left\langle u_{f},\phi_{j}\right\rangle $ and $\left\langle v_{f},\phi_{j}\right\rangle $ due to the low level of frequency restricted by the stabilized operator. It then enables us to obtain the reconstructions $U_{f}^{\varepsilon,j},V_{f}^{\varepsilon,j}$.

We now take into account the numerical regularized solution by using the backward Euler method. In this regard, a uniform grid of mesh-points $(x,t) = (x_m,t_k)$ is used. Here $x_m = a + m\Delta x$ and $t_k = k\Delta t$ for $0\le m\le M$ and $0\le k\le K$ where $\Delta x = \frac{b-a}{M}$,  $\Delta t = \frac{T}{K}$ and $M,K\in\mathbb{N}^{*}$. Henceforth, we find $u_{m,k}^{\varepsilon,n}\approx u^{\varepsilon,n}(x=x_m,t=t_k)$ and $v_{m,k}^{\varepsilon,n}\approx v^{\varepsilon,n}(x=x_m,t=t_k)$ by solving the following fully discretized system:
\begin{align*}
	\frac{u_{m,k+1}^{\varepsilon,n}-u_{m,k}^{\varepsilon,n}}{\Delta t}+\overline{\mathcal{D}}_{1}\left(u_{k+1}^{\varepsilon,n}\right) & \frac{u_{m+1,k}^{\varepsilon,n}-2u_{m,k}^{\varepsilon,n}+u_{m-1,k}^{\varepsilon,n}}{\Delta x^{2}}\\
	& =F_{\ell^{\varepsilon}}\left(u_{m,k+1}^{\varepsilon,n};v_{m,k+1}^{\varepsilon,n}\right)+\mathbf{P}_{\varepsilon}^{\beta}u_{m,k+1}^{\varepsilon,n},\\
	\frac{v_{m,k+1}^{\varepsilon,n}-v_{m,k}^{\varepsilon,n}}{\Delta t}+\overline{\mathcal{D}}_{2}\left(v_{k+1}^{\varepsilon,n}\right) & \frac{v_{m+1,k}^{\varepsilon,n}-2v_{m,k}^{\varepsilon,n}+v_{m-1,k}^{\varepsilon,n}}{\Delta x^{2}}\\
	& =G_{\ell^{\varepsilon}}\left(u_{m,k+1}^{\varepsilon,n};v_{m,k+1}^{\varepsilon,n}\right)+\mathbf{P}_{\varepsilon}^{\beta}v_{m,k+1}^{\varepsilon,n},
\end{align*}
for $1\le m \le M-1$ and $0\le k \le K-1$. Accordingly, we endow this system with the Dirichlet conditions $u^{\varepsilon,n}_{0,k}=u^{\varepsilon,n}_{M,k}=v^{\varepsilon,n}_{0,k}=v^{\varepsilon,n}_{M,k}=0$ and the terminal conditions
$u^{\varepsilon,n}_{m,K}=U^{\varepsilon,n}_{f}(x=x_m)$, $v^{\varepsilon,n}_{m,K}=V^{\varepsilon,n}_{f}(x=x_m)$. Note that our objective here is to solve this system backwards in time and thus, $k$ runs from $K$ to $0$ in the iterations. Then, for each $k=0,\ldots,K-1$ we arrive at a linear algebraic system with a Toeplitz-like matrix, as follows:
\begin{align}\label{eq:algebra}
\mathbb{A}\left(u_{k+1}^{\varepsilon,n},v_{k+1}^{\varepsilon,n}\right)\mathbb{W}_{k}^{\varepsilon,n}=\mathbb{F}_{\varepsilon}\left(u_{k+1}^{\varepsilon,n};v_{k+1}^{\varepsilon,n}\right),
\end{align}
where the solution $\mathbb{W}^{\varepsilon,n}_{k}\in \mathbb{R}^{2(M-1)}$, the block matrix $\mathbb{A}\in\mathbb{R}^{(2M-2)\times(2M-2)}$ and the function $\mathbb{F}_{\varepsilon}\in \mathbb{R}^{2(M-1)}$ are defined by
\[
\mathbb{W}_{k}^{\varepsilon,n}=\begin{bmatrix}u_{1,k}^{\varepsilon,n} & u_{2,k}^{\varepsilon,n} & \cdots & u_{M-2,k}^{\varepsilon,n} & u_{M-1,k}^{\varepsilon,n} & v_{1,k}^{\varepsilon,n} & v_{2,k}^{\varepsilon,n} & \cdots & v_{M-2,k}^{\varepsilon,n} & v_{M-1,k}^{\varepsilon,n}\end{bmatrix}^{\text{T}},
\]
\[
\mathbb{F}_{\varepsilon}\left(u_{k+1}^{\varepsilon,n};v_{k+1}^{\varepsilon,n}\right)=\begin{bmatrix}u_{1,k+1}^{\varepsilon,n}-\Delta tF_{\ell^{\varepsilon}}\left(u_{1,k+1}^{\varepsilon,n};v_{1,k+1}^{\varepsilon,n}\right)-\Delta t\mathbf{P}_{\varepsilon}^{\beta}u_{1,k+1}^{\varepsilon,n}\\
u_{2,k+1}^{\varepsilon,n}-\Delta tF_{\ell^{\varepsilon}}\left(u_{2,k+1}^{\varepsilon,n};v_{2,k+1}^{\varepsilon,n}\right)-\Delta t\mathbf{P}_{\varepsilon}^{\beta}u_{2,k+1}^{\varepsilon,n}\\
\vdots\\
u_{M-2,k+1}^{\varepsilon,n}-\Delta tF_{\ell^{\varepsilon}}\left(u_{M-2,k+1}^{\varepsilon,n};v_{M-2,k+1}^{\varepsilon,n}\right)-\Delta t\mathbf{P}_{\varepsilon}^{\beta}u_{M-2,k+1}^{\varepsilon,n}\\
u_{M-1,k+1}^{\varepsilon,n}-\Delta tF_{\ell^{\varepsilon}}\left(u_{M-1,k+1}^{\varepsilon,n};v_{M-1,k+1}^{\varepsilon,n}\right)-\Delta t\mathbf{P}_{\varepsilon}^{\beta}u_{M-1,k+1}^{\varepsilon,n}\\
v_{1,k+1}^{\varepsilon,n}-\Delta tG_{\ell^{\varepsilon}}\left(u_{1,k+1}^{\varepsilon,n};v_{1,k+1}^{\varepsilon,n}\right)-\Delta t\mathbf{P}_{\varepsilon}^{\beta}v_{1,k+1}^{\varepsilon,n}\\
v_{2,k+1}^{\varepsilon,n}-\Delta tG_{\ell^{\varepsilon}}\left(u_{2,k+1}^{\varepsilon,n};v_{2,k+1}^{\varepsilon,n}\right)-\Delta t\mathbf{P}_{\varepsilon}^{\beta}v_{2,k+1}^{\varepsilon,n}\\
\vdots\\
v_{M-2,k+1}^{\varepsilon,n}-\Delta tG_{\ell^{\varepsilon}}\left(u_{M-2,k+1}^{\varepsilon,n};v_{M-2,k+1}^{\varepsilon,n}\right)-\Delta t\mathbf{P}_{\varepsilon}^{\beta}v_{M-2,k+1}^{\varepsilon,n}\\
v_{M-1,k+1}^{\varepsilon,n}-\Delta tG_{\ell^{\varepsilon}}\left(u_{M-1,k+1}^{\varepsilon,n};v_{M-1,k+1}^{\varepsilon,n}\right)-\Delta t\mathbf{P}_{\varepsilon}^{\beta}v_{M-1,k+1}^{\varepsilon,n}
\end{bmatrix},
\]
\[
\mathbb{A}\left(u_{k+1}^{\varepsilon,n},v_{k+1}^{\varepsilon,n}\right)=\begin{bmatrix}\mathbb{D}_{1}^{k+1} & O\\
O & \mathbb{D}_{2}^{k+1}
\end{bmatrix},
\]
with $\mathbb{D}_{1}^{k+1},\mathbb{D}_{2}^{k+1}\in \mathbb{R}^{K-1}$ taken by
\[
\mathbb{D}_{1}^{k+1}=\begin{bmatrix}1+2\overline{\mathcal{D}}_{1}^{k+1}\bar{\alpha} & -\overline{\mathcal{D}}_{1}^{k+1}\bar{\alpha} & \cdots & \cdots & 0\\
-\overline{\mathcal{D}}_{1}^{k+1}\bar{\alpha} & 1+2\overline{\mathcal{D}}_{1}^{k+1}\bar{\alpha} & -\overline{\mathcal{D}}_{1}^{k+1}\bar{\alpha} &  & \vdots\\
\vdots &  &  &  & \vdots\\
\vdots &  & -\overline{\mathcal{D}}_{1}^{k+1}\bar{\alpha} & 1+2\overline{\mathcal{D}}_{1}^{k+1}\bar{\alpha} & -\overline{\mathcal{D}}_{1}^{k+1}\bar{\alpha}\\
0 & \cdots & \cdots & -\overline{\mathcal{D}}_{1}^{k+1}\bar{\alpha} & 1+2\overline{\mathcal{D}}_{1}^{k+1}\bar{\alpha}
\end{bmatrix},
\]
\[
\mathbb{D}_{2}^{k+1}=\begin{bmatrix}1+2\overline{\mathcal{D}}_{2}^{k+1}\bar{\alpha} & -\overline{\mathcal{D}}_{2}^{k+1}\bar{\alpha} & \cdots & \cdots & 0\\
-\overline{\mathcal{D}}_{2}^{k+1}\bar{\alpha} & 1+2\overline{\mathcal{D}}_{2}^{k+1}\bar{\alpha} & -\overline{\mathcal{D}}_{2}^{k+1}\bar{\alpha} &  & \vdots\\
\vdots &  &  &  & \vdots\\
\vdots &  & -\overline{\mathcal{D}}_{2}^{k+1}\bar{\alpha} & 1+2\overline{\mathcal{D}}_{2}^{k+1}\bar{\alpha} & -\overline{\mathcal{D}}_{2}^{k+1}\bar{\alpha}\\
0 & \cdots & \cdots & -\overline{\mathcal{D}}_{2}^{k+1}\bar{\alpha} & 1+2\overline{\mathcal{D}}_{2}^{k+1}\bar{\alpha}
\end{bmatrix},
\]
for $\overline{\mathcal{D}}_{1}^{k+1}=\overline{\mathcal{D}}_{1}\left(u_{k+1}^{\varepsilon,n}\right)$,
$\overline{\mathcal{D}}_{2}^{k+1}=\overline{\mathcal{D}}_{2}\left(v_{k+1}^{\varepsilon,n}\right)$ and $\bar{\alpha}=\frac{\Delta t}{\Delta x^2}$.

We remark that for each iteration step $k$, \eqref{eq:algebra} has a solution because the block matrix $\mathbb{A}$ is symmetric positive definite. In this illustration,
we only seek the approximation of $u^{\text{ex}}(x_m,t=t^{\varepsilon})$, but usually $t^{\varepsilon}$ does not belong to the underlying mesh-points. We then rely on the nearest point of this $t^{\varepsilon}$, which means that for each $\varepsilon>0$ the error, denoted by $E(t_k)$, is computed for $t_k$ such that $t_k\ge t^{\varepsilon}$ and $t_{k-1}< t^{\varepsilon}$. Given the number of samples $\mathcal{M}\in \mathbb{N}^{*}$ and choosing $T=1$, the averaged error $E(t_k)$ reads as 
\[
E(t_k) = \frac{1}{\mathcal{M}}  \frac{1}{M} \sum\limits_{m_s=1}^{\mathcal M} \sum\limits_{m=1}^M \left | u_{m,k,m_s}^{\varepsilon,n}- u^{\text{ex}}(x_m,0)\right|^2,
\]
where $u^{\text{ex}}$ denotes the manufactured solution that we will choose in the numerical tests below. Note that for simplicity, we eventually omit the index $m_s$ in the regularized solution in the above numerical setting. The value of $t^{\varepsilon}$ can be approximately computed by the equation $
\log(t^{\varepsilon})= -t^{\varepsilon} \log\left(\varepsilon^{\frac{1}{T}\max\left\{\theta -1, -4\theta p \right\}}\right)
$ (cf. \cref{sec:proof2}).



\subsection{Test 1} We begin the numerical verification by the most widely used Lotka--Volterra model for tumor growth kinetics where in \eqref{system1} one has $\mathcal{D}_{i}$ are positive constants and the corresponding Lotka--Volterra functions are given by
\begin{align*}
	F = r_{u}u\left(1 - \frac{u}{K_{u}} - \frac{b_{u}v}{K_{u}}\right) + F_1 (x,t),\quad  G = r_{v}v\left(1 - \frac{v}{K_{v}} - \frac{b_{v}u}{K_{v}}\right) + F_2 (x,t).
\end{align*}
Here, the involved parameters in $F,G$ are positive constants and cf. \cite{Gatenby2002}, they possess the following biological meaning: $b_{u}$ is the negative effects of tumor on normal cells, $b_{v}$ is  a lumped, phenomenological term presented by, e.g., the immune response; $r_{u},r_{v}$ are, respectively, maximum growth rates of normal and tumor cells; $K_{u}, K_{v}$ denote their maximal densities.

Given $\mathcal{D}_{1}=0.5$ and $\mathcal{D}_{2}=1$ for slow diffusions, which then lead to $\overline{M}=1.1$ and $M_1=1$, we choose $r_u = r_v = K_u = K_v = 1$ and $b_u = 2, b_v = 0.5$, usually describing the slight domination of tumor cells. The manufactured solutions are given by $u^{\text{ex}}=e^{-t}\sin(x)$ and $v^{\text{ex}}=t^2 x(\pi-x)$, ensuring the positivity of solutions. This way we need to take
\begin{align*}
	F_{1}\left(x,t\right) & =e^{-t}\sin\left(x\right)\left[e^{-t}\sin\left(x\right)+2t^{2}x\left(\pi-x\right)-\frac{3}{2}\right],\\
	F_{2}\left(x,t\right) & =2t^{2}+tx\left(\pi-x\right)\left[2-t+t^{3}x\left(\pi-x\right)+\frac{t}{2}e^{-t}\sin\left(x\right)\right].
\end{align*}

By the choice of $F$ and $G$, we find that $L(\ell^{\varepsilon}) = 1 + 3\ell^{\varepsilon}$ and hence, for every time step $\ell^{\varepsilon}$ can be controlled by
\[
\ell^{\varepsilon}\left(t\right)=\frac{1}{3}\left(\frac{1}{16\left(T-t\right)}\log\left(\log^{\kappa\left(t\right)}\left(\varepsilon^{\frac{1}{C_{1}T}\max\left\{ \theta-1,-4\theta p\right\} }\right)\right)-1\right),
\]
where we choose $\kappa(t)=C_1 t>0$. Thus, we can work with the cut-off functions $F_{\ell^{\varepsilon}}$ and $G_{\ell^{\varepsilon}}$ in this circumstance.

\subsection{Test 2}
In the second example, we take into account the nonlocal problem as in \cite{Chipot1997}, where $F,G$ are independent of the solutions, i.e. $F(x,t;u;v)=F(x,t)$ and $G(x,t;u;v)=G(x,t)$. Then the diffusion terms are particularly expressed by
\[
D_i(u)(t) = d_i + \int\limits_{\Omega, u\in [0,1]} u(1-u)dx.
\]

Basically, this context can be understood that the tumor cells are on the way to dominate the normal cells, since the evolution system is now uncoupled and partly driven by the total population of each species. In this test, we choose $u^{\text{ex}} = 0.25\log(3+t)\sin^2(x)$ and $v^{\text{ex}}=t\sin(x)$ with $d_1 = 0.01$ and $d_2 = 0.05$. Then, it follows that $\overline{M}= 4, M_1 = 3.5$ and
\begin{align*}
	F\left(x,t\right) & =\frac{1}{4\left(3+t\right)}\sin^{2}\left(x\right)
	\\ &
	-\frac{1}{8}\left(\frac{1}{25}+\frac{\pi}{2}\log\left(3+t\right)\left(1-\frac{3}{16}\log\left(3+t\right)\right)\right)\log\left(3+t\right)\cos\left(2x\right),\\
	G\left(x,t\right) & =\left(1-t\left(\frac{t\left(\pi t-4\right)}{2}-\frac{1}{20}\right)\right)\sin\left(x\right).
\end{align*}

Note that due to the nonlocal diffusion, we compute $\overline{\mathcal{D}}_{i}^{k+1}$ ($i=1,2$) involved in the algebraic system \eqref{eq:algebra} by employing the midpoint rule. It yields
\[
\overline{\mathcal{D}}_{i}^{k+1}=\overline{M}-\mathcal{D}_{i}\left(u_{k+1}^{\varepsilon,n}\right)\approx
\overline{M}-d_{i}-\frac{b-a}{M}\sum_{m=1}^{M-1}\left.\tau\left(u_{m,k+1}^{\varepsilon,n}\right)\right|_{u_{m,k+1}^{\varepsilon,n}\in\left[0,1\right]},
\]
where $\tau\left(u\right)=\left|u\left(1-u\right)\right|$. It is worth noting that the choice of the numerical integration for the diffusion term has to be careful. Although the efficient Gauss--Legendre quadrature and its variants seem applicable by choosing suitable nodes, in practice one cannot get measurements exactly at those abscissae. One can also take into account the Newton--Cotes method for the uniform grid in this implementation, but this approach is usually unstable when $M$ is large. Therefore, one can see the impediment in handling regularization methods in numerics, especially when there are numerous nonlinear factors involved in the model.

\subsection{Numerical results} We take $p = 1$ and $\theta = 0.3$ in these two tests. The choice of $\theta$ indicates that we aim to see the convergence for a not very large number of steps of discrete observations $n$. Meanwhile, the value of $
\varepsilon$ varies from $10^{-3}$ to $10^{-5}$.

\begin{figure}[htbp]
	\hspace*{-1.7cm}\includegraphics[width=1.12\textwidth]{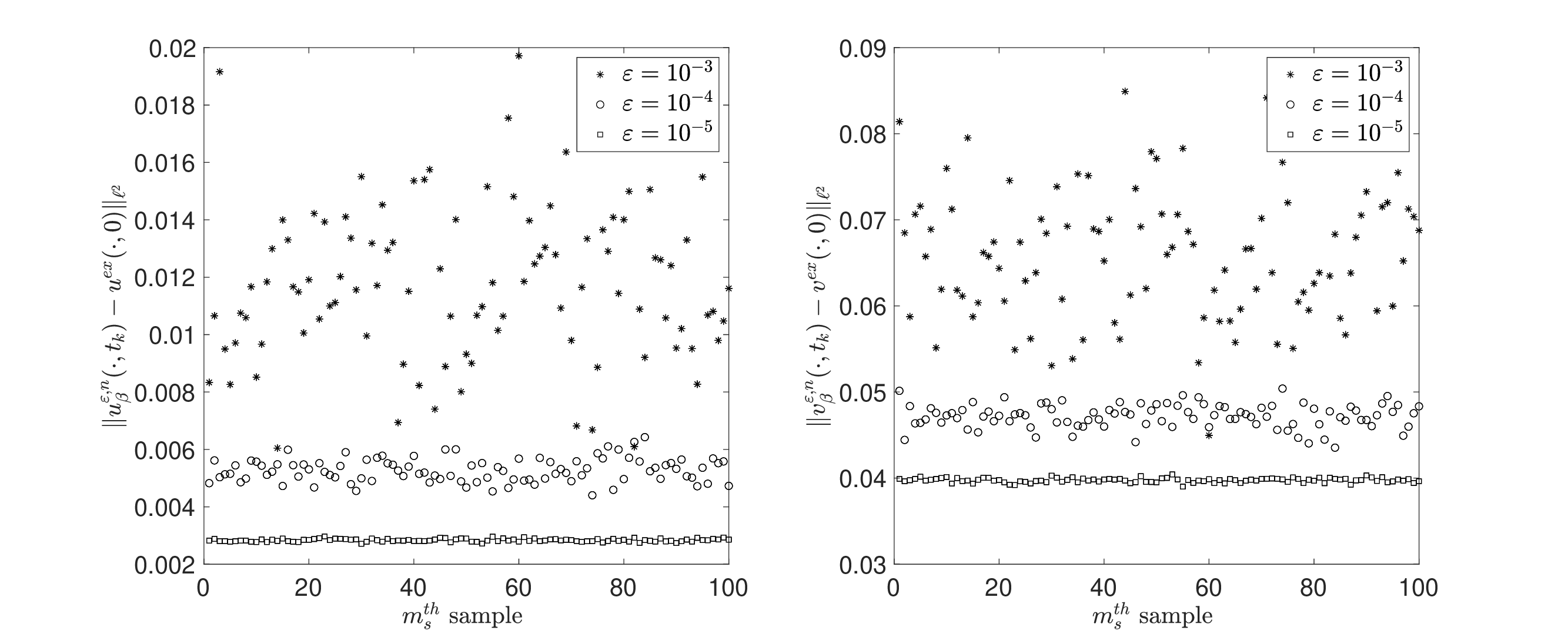}
	\caption{Numerical errors between the regularized solution and the true solution in Test 1 with 100 samples and $\varepsilon\in \left\{10^{-3},10^{-4},10^{-5}\right\}$.}
	\label{fig:1}
\end{figure}

\begin{figure}[htbp]
	\hspace*{-1.7cm}\includegraphics[width=1.12\textwidth]{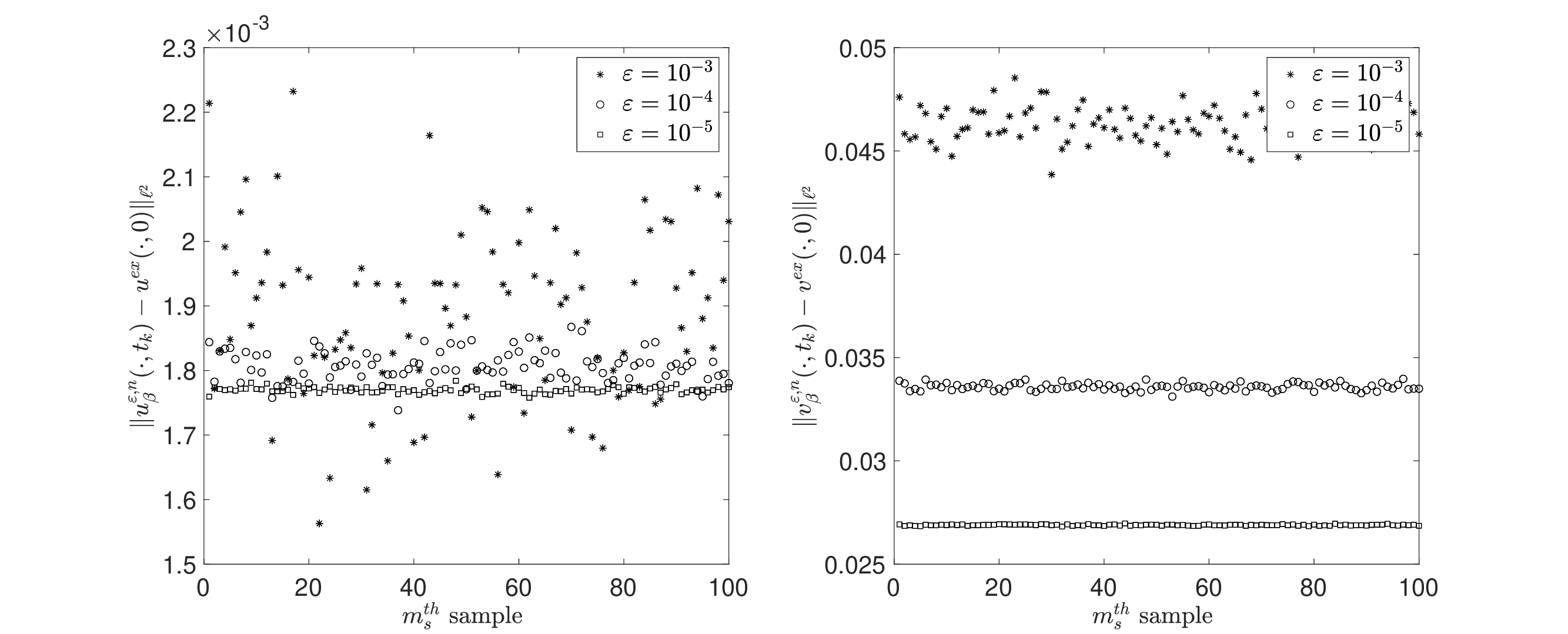}
	\caption{Numerical errors between the regularized solution and the true solution in Test 2 with 100 samples and $\varepsilon\in \left\{10^{-3},10^{-4},10^{-5}\right\}$.}
	\label{fig:2}
\end{figure}

Our convergence results are concluded in \cref{fig:1,fig:2}, where we have computed the numerical errors of two tests with $\mathcal{M}=100$ samples (the number of samples is implemented in the horizontal line of each figure) and $M=15$, $K=100$. Having these errors allows us to obtain the averaged errors corresponding to each solution, which afterward confirms our theoretical expectation. Particularly, in Test 1 the averaged error for $u$ decreases from $0.01224$ (for $\varepsilon = 10^{-3}$) to $0.00284$ (for $\varepsilon = 10^{-5}$), while for $v$ it slowly goes down to $0.03971$ from $0.06545$. Similarly, in Test 2 the averaged error for $u$ is from $0.00191$ (for $\varepsilon = 10^{-3}$) to $0.00177$ (for $\varepsilon = 10^{-5}$), while for $v$ it decreases from $0.04640$ to $0.02689$. Henceforth, one can see that our approximation is acceptable from $\varepsilon = 10^{-3}$, i.e. this time the regularized solution is close to the true solution.

\section{Conclusions}
\label{sec:conclusions}

We have improved and adapted the recently developed QR method in \cite{Nguyen2019} to regularization of a terminal-boundary value parabolic problem with white Gaussian noise. Although this noise process is very standard in the stochastic setting, this attempt concretely consolidates the mathematical quality of the reconstruction method we have established so far, where both the well-posedness of the regularized problem and the convergence analysis are successfully obtained. Besides, it again confirms the conventional logarithmic rate that we usually meet in regularization of parabolic problems. This time we introduce a new type of perturbation that allows us to compute the stabilized operator in the numerical verification of the scheme.

Although several questions remain open in \cite{Nguyen2019}, we would like to add one more interesting problem. In view of the localization of the brain tumor source in this paper, one can address the model in a domain with moving boundaries as investigated in, e.g., \cite{Almeida2018}. This way may be helpful in brain pathologies, especially in improving the Magnetic Resonance Imaging method as one wants to anatomically see the motion of the brain (see, e.g., \cite{Holdsworth2016}).



\section*{Acknowledgments}
V.A.K thanks Thi Kim Thoa Thieu (Gran Sasso Science Institute, Italy and Karlstad University, Sweden) for fruitful discussions during the completion of the manuscript. V.A.K thanks Prof. Loc Nguyen (North Carolina, USA) for the hospitality he provided.

\bibliographystyle{siamplain}
\bibliography{references}

\end{document}